\documentclass[12pt,a4paper]{article}
\usepackage{theorem}
\usepackage{epsfig}
\usepackage{amsmath}
\usepackage{amssymb}
\usepackage{amsfonts}
\usepackage{pst-node}
\usepackage{html}

\newtheorem{thm}{Theorem}

\newtheorem{lemma}{Lemma}
\newtheorem{note}{Note}
\newtheorem{cor}{Corollary}
\theorembodyfont{\rmfamily}

\newenvironment{proof}
{\begin{rm}\par\smallskip\noindent{\bf Proof.}\quad}{\QED\end{rm}}
\newenvironment{ack}
{\begin{rm}\par\bigskip\noindent{\bf Acknowledgements.}\quad\em}{\end{rm}}
\newcommand{\conv}{\mathop{\mathrm{conv}}\nolimits}
\def\rel{\mathop{\mathrm{relint}}\nolimits}
\def\cl{\mathop{\mathrm{cl}}\nolimits}

\title{
On f-vectors of Minkowski additions of convex polytopes%
\footnote{Supported by the Swiss National Science Foundation Project 200021-105202,
``Polytopes, Matroids and Polynomial Systems''.}
}
\date{October 21, 2005\\Revised October 24, 2006}
\author{Komei Fukuda%
\footnote{Also affiliated with Institute for Operations Research
and Institute of Theoretical Computer Science
ETH Zentrum, Zurich, Switzerland.}
\\
Mathematics Institute\\
EPFL, Lausanne\\Switzerland\\
komei.fukuda@epfl.ch\\
\and 
Christophe Weibel\\
Mathematics Institute\\
EPFL, Lausanne\\Switzerland\\
christophe.weibel@epfl.ch
}

\def\NN{\mathcal{N}}
\def\RRd{\mathbb{R}^d}
\def\BBox{\rule{2mm}{3mm}}
\def\QED{\hfill$\BBox$}

\begin{document}
\maketitle

\begin{abstract}

The objective of this paper is to present two types of results
on Minkowski sums of convex polytopes.
The first is about a special class of polytopes we call
\emph{perfectly centered} and the combinatorial properties
of the Minkowski sum with their own dual. In particular, we 
have a characterization of face lattice of the sum  in terms of
the face lattice of a given perfectly centered polytope. Exact face
counting formulas are then obtained for perfectly centered simplices
and hypercubes. The second type of results concerns tight \emph{upper bounds}
for the f-vectors of Minkowski sums of several polytopes.
\end{abstract}
\section{Introduction}

Minkowski sums of polytopes in $\RRd$ naturally arise in many domains,
ranging from mechanical engineering
\cite{Petit04} to algebra \cite{Sturmfels96,Sturmfels98}.
These applications have triggered recent algorithmic advances 
\cite{Fukuda04} and an efficient implementation \cite{w-msv-05}.  
Despite the new developments, we are still very far from
understanding the combinatorial structure 
(i.e. the face lattice) of a Minkowski sum of several polytopes.
In particular, it is in general 
difficult to estimate the number of $k$-dimensional faces
($k$-faces) of the result for each $0\le k\le d-1$, 
even if we know the face lattices of the summands.
One special case of the problem which is relatively well
understood is when the summands are $m$ line segments
in $\RRd$.  The resulting sum is known as a {\em zonotope\/} given by $m$ generators, see e.g. \cite[Lecture 7]{Ziegler95}.
The goal of the paper is to study the problem, first in the particular
case of a certain class of polytopes summed with their own dual, and
then to prove tight upper bounds on the number of $k$-faces.

We call a polytope {\em centered\/} if
 it contains the origin in its relative interior.
Nesterov \cite{Nesterov04} has recently proved 
that the sum of a centered full-dimensional polytope (and more generally
a centered full-dimensional compact convex body) with its dual, if properly scaled,
gives a set whose asphericity is at most the square root of
that of the initial polytope. The \emph{asphericity} of a set is here
defined as the ratio of the diameter of its smallest enclosing
ball to that of its largest enclosed ball. Thus, summing a polytope
with its own dual has a strong \emph{rounding effect}.  For this reason, 
the Minkowski sum $P + \alpha P^*$, will be called
a \emph{Nesterov rounding} of a polytope $P$ for any positive
scalar $\alpha$. 

Of special interest is the combinatorial
aspect of Nesterov rounding. The first observation  is
that the combinatorial structure 
of $P + \alpha P^*$ does not depend on $\alpha$. Thus 
we can set the scaling factor to be $1$ without loss of
generality.  We shall say ``the'' Nesterov rounding instead
of ``a'' Nesterov rounding to mean the class of all
Nesterov roundings with the unique combinatorial type.
While the scaling factor is irrelevant for our study,
the position of the origin in $P$ does affect the combinatorial
structure of the Nesterov rounding.  In other words,
the combinatorial structure of the Nesterov rounding of $P$
is not uniquely determined by that of the polytope. 
However, it is the case when the polytope has the
\emph{perfectly centered\/} property, by which we mean
that every nonempty face intersects with
its outer normal cone, see Section~\ref{pcp} for the formal
definition.

The first result characterizes the face lattice of the Nesterov
rounding of a perfectly centered polytope.
For this, we use the natural bijection between the faces of $P$
and those of the dual: $F^D$ denotes the dual face associated
with a face $F$ of $P$.  We call a face $F$ of a polytope $P$ {\em trivial\/}
if it is either the empty set $\emptyset$ or the polytope $P$ itself.  In particular,
the trivial faces are dual to each other: $P^D =\emptyset$.

\begin{thm}\label{mainthm}
Let $P$ be a perfectly centered polytope.
A subset $H$ of $P+P^*$ is a nontrivial face of  $P+P^*$ if 
and only if  $H=G+F^D$ for some ordered nontrivial faces 
$G\subseteq F$ of $P$.
\end{thm}

This theorem can be considered as a natural
 extension of a theorem in
\cite{Broadie85} which was restricted to 
the facets of the sum.

As a corollary, we obtain face counting formulas for 
perfectly centered simplices (Theorem \ref{thm:simplexrounding}) 
and hypercubes (Theorem \ref{thm:cuberounding}).

Any face of a Minkowski sum of polytopes can be \emph{decomposed}
uniquely into a sum of faces of the summands. We will say that the
decomposition is \emph{exact} when the dimension of the sum is
equal to the sum of the dimensions of the summands. When all facets
have an exact decomposition, we will say the summands are
\emph{relatively in general position}.

This provides us with a trivial upper bound
for the number of faces, i.e. the number of possible distinct
decompositions. As usual, we denote by $f_k(P)$ the number of $k$-faces
of a $d$-polytope $P$. For each $k=0,\ldots,d-1$ and $n\geq 1$,
the number of $k$-faces of $P_1+\cdots+P_n$ is bounded by:
$$
f_k(P_1+\cdots+P_n) \leq
\sum_{
\begin{array}{c}
1\leq s_i \leq f_0(P_i)\\
s_1+\cdots+s_n=k+n
\end{array}}
\prod_{i=1}^n
\left(\begin{array}{c}
f_0(P_i)\\
s_i
\end{array}\right),
$$
where $s_i$'s are integral.

The next theorem shows that this bound can be
achieved in some cases.

\begin{thm}\label{thmvertices}
In dimension $d\geq 3$, it is possible to choose $d-1$ polytopes so that
the trivial upper bound for the number of vertices is attained.
\end{thm}

Tight upper bounds on the number of facets appear to be 
harder to obtain in general. However, the following result
on 3-polytopes holds.

\begin{thm}\label{thmfaces}
Let $P_1,\ldots,P_n$ be 3-dimensional polytopes relatively in general
position, and $P$ their sum. Then the following equations hold:
$$
2 f_2(P)-f_1(P) = \sum_{i=1}^n (2 f_2(P_i)- f_1(P_i)) ,
$$
$$
f_2(P)-f_0(P)+2 = \sum_{i=1}^n (f_2(P_i)-  f_0(P_i)+2) ,
$$
$$
f_1(P)-2f_0(P)+4 = \sum_{i=1}^n (f_1(P_i)-2 f_0(P_i)+4 ).
$$
\end{thm}
As a corollary (Corollary~\ref{cor:facets}), we obtain tight upper bounds for the number of facets (and edges) of the sum of two $3$-polytopes.

Furthermore, by using the fact that the Nesterov rounding
of a perfectly centered polytope is again perfectly centered
(Theorem~\ref{thm:pcpc}), 
we also analyze the asymptotic behavior of repeated Nesterov
roundings in dimension 3 (Theorem \ref{thm:repeatedNR}) that
in fact shows a combinatorial rounding effect: the ratio of the number of vertices
over that of facets approaches $1$.

When the dimension $d$ is large enough relative to the number of polytopes,
it is possible for faces of lower dimensions to attain the trivial upper bounds.

\begin{thm} \label{thmcyc}
In dimension $d\geq 4$,
it is possible to choose $n\leq\lfloor\frac{d}{2}\rfloor$
polytopes $P_1$, $P_2$, $\ldots$, $P_n$ so that the trivial upper bound for 
the number of $k$-faces
of $P_1+ \cdots + P_n$ is attained for all $0\le k \le \lfloor\frac{d}{2}\rfloor-n$.
\end{thm}

Throughout this paper, we assume that the reader is familiar with
the basic results on convex polytopes.
For a general introduction to polytopes,
please refer to \cite{Grunbaum67,Ziegler95}.

\section{Perfectly Centered Polytopes}\label{pcp}

We will assume in this section that all polytopes are full-dimensional.
A polytope $P$ is said to be \emph{centered} if its relative interior $\rel(P)$
contains the origin. For any centered polytope $P$, its \emph{dual},
denoted by $P^*$, is defined by
$$
P^* = \{x\;:\;\langle x,y\rangle \leq 1,\;\forall y\in P\}.
$$
For each face of $P$, we define the
\emph{associated dual face} as
$$
F^D = \left\{x \;|\; x\in P^*\;:\;\langle x, f\rangle =1,\forall f\in F\right\} .
$$
We will define now two notions central to Minkowski sums.

For a polytope $ P$ in $\RRd$ and any vector $c\in\RRd$, we denote
by $S(P;c)$
the \emph{set of maximizers} of the linear function $\langle\cdot,c\rangle$:
$$
S(P;c)\;=\;\left\{x\in P \;|\; \langle x, c\rangle =\max_{y\in P}\langle y, c\rangle \right\} .
$$
For any face $F$ of $P$, the \emph{outer normal cone} of $P$ at $F$, denoted by $\NN(F;P)$, 
is the set
of vectors $c$ such that $F=S(P;c)$.
Normal cones are relatively open. Also, if $F$ and $G$ are nonempty faces of a polytope $P$,
$$
G\subseteq F \Leftrightarrow
\cl(\NN(F;P)) \subseteq \cl(\NN(G;P)),
$$
where $\cl(S)$ denotes the \emph{topological closure} of a set $S$.

\begin{lemma}\label{coneface}
Let $P$ be a centered polytope. For a face $F$ of $P$, $F^D$ is
a face of $P^*$. Furthermore, if $F$ is nontrivial,
$\NN(F;P)$ is the cone generated by the points in
the relative interior of the dual face $F^D$.  Namely,
\[
\NN(F;P)=\{\lambda x\;:\;\lambda>0,\;x\in \rel(F^D)\}.
\]
Consequently,
$$
\cl(\NN(F;P))=\{\lambda x\;:\;\lambda\geq 0,\;x\in F^D\}.
$$
\begin{proof}
The proof is straightforward, and left to the reader.
\end{proof}
\end{lemma}

\begin{cor}
Let $F$ be a nontrivial face of a polytope $P$.
The affine spaces spanned by $F$ and $F^D$ are orthogonal to each other,
meaning, the linear subspaces obtained from the affine spaces by
translations are orthogonal.
\end{cor}

The study of Nesterov rounding of polytopes has led to a new class
of polytopes, that we introduce now. A polytope is called
\emph{perfectly centered} if
\[
\rel(F)\cap \NN(F;P) \neq \emptyset
\text{ for any nonempty face } F \text{ of } P.
\]
Observe that if the intersection is nonempty, then it consists of a single point,
since a face is orthogonal to its normal cone.
\begin{figure}[htb]
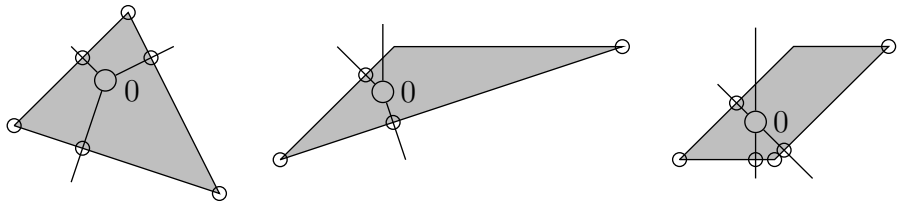

$$
\begin{array}{ccc}
\psset{xunit=0.6cm,yunit=0.6cm,shortput=nab,linewidth=0.5pt,arrowsize=2pt 3,labelsep=2cm}
\pspicture[0.5](-2.5,-3)(3,2)
\psline[fillstyle=solid,fillcolor=lightgray,linetype=1](-2,-1)(2.5,-2.5)(0.5,1.5)(-2,-1)
\cnode(-2,-1){0.1}{G}
\cnode(2.5,-2.5){0.1}{H}
\cnode(0.5,1.5){0.1}{I}
\cnode(0,0){0.15}{O}
\cnode(-0.5,-1.5){0.1}{A}
\pnode(-0.75,-2.25){B}
\cnode(1,0.5){0.1}{C}
\pnode(1.5,0.75){D}
\cnode(-0.5,0.5){0.1}{E}
\pnode(-0.75,0.75){F}
\ncline{O}{B}
\ncline{O}{D}
\ncline{O}{F}
\nput{d}{0,0}{$\;\;\;\;\;\;0$}
\endpspicture
&
\psset{xunit=0.6cm,yunit=0.6cm,shortput=nab,linewidth=0.5pt,arrowsize=2pt 3,labelsep=2pt}
\pspicture[0.5](-2.5,-2)(5.5,1.5)
\psline[fillstyle=solid,fillcolor=lightgray,linetype=1](-2.25,-1.5)(5.25,1)(0.25,1)(-2.25,-1.5)
\cnode(-2.25,-1.5){0.1}{G}
\cnode(5.25,1){0.1}{H}
\cnode(0,0){0.15}{O}
\cnode(0.225,-0.675){0.1}{A}
\pnode(0.5,-1.5){B}
\pnode(0,1.5){D}
\cnode(-0.375,0.375){0.1}{E}
\pnode(-1,1){F}
\ncline{O}{B}
\ncline{O}{D}
\ncline{O}{F}
\nput{r}{0,0}{$\;\;0$}
\endpspicture
&
\psset{xunit=0.5cm,yunit=0.5cm,shortput=nab,linewidth=0.5pt,arrowsize=2pt 3,labelsep=2pt}
\pspicture[0.5](-2.5,-1.5)(4,2.5)
\psline[fillstyle=solid,fillcolor=lightgray,linetype=1](-2,-1)(0.5,-1)(3.5,2)(1,2)(-2,-1)
\cnode(-2,-1){0.1}{I}
\cnode(0.5,-1){0.1}{J}
\cnode(3.5,2){0.1}{K}
\cnode(0,0){0.15}{O}
\cnode(0,-1){0.1}{A}
\pnode(0,-1.5){B}
\cnode(0.75,-0.75){0.1}{C}
\pnode(1.5,-1.5){D}
\pnode(0,2.5){F}
\cnode(-0.5,0.5){0.1}{G}
\pnode(-1,1){H}
\ncline{O}{B}
\ncline{O}{D}
\ncline{O}{F}
\ncline{O}{H}
\nput{r}{0,0}{$\;\;0$}
\endpspicture
\end{array}
$$
\caption{\label{fig:imperf}A perfectly centered 
and two non-perfectly-centered polytopes}
\end{figure}

For instance, the polytope on the left in Figure~\ref{fig:imperf} is perfectly centered,
and the two others are not. 
The one in the center can be made
perfectly centered by moving the origin, but the one on the right cannot be.
Note that the perfectly centered property was
previously studied in  \cite{Broadie85}
where it was called the {\em projection condition}. 
Advantages of using the term ``perfectly centered'' 
over the old term become evident when we state theorems such as
Corollary \ref{pcduality} and Theorem \ref{thm:pcpc} below.
\begin{lemma}\label{lem:coneinter}
A polytope $P$ is perfectly centered if and only if $P$ is centered and
$\NN(F;P)\cap\NN(F^D;P^*)\neq\emptyset$,
for every nontrivial face $F$ of $P$.
\begin{proof}
By Lemma~\ref{coneface}, for any nontrivial face $F$,
$$
\NN(F^D;P^*)=\{\lambda x\;:\;\lambda>0,\;x\in \rel(F)\}.
$$
Thus, for every nontrivial face $F$ of a polytope $P$, the relations
$\rel(F)\cap\NN(F;P)\neq\emptyset$ and
$\NN(F^D;P^*) \cap\NN(F;P)\neq\emptyset$ are equivalent.
Since $\NN(P;P)=\{0\}$, two statements
$\rel(P)\cap\NN(P;P)\neq\emptyset$ and  $0\in \rel(P)$ are also
equivalent.
\end{proof}
\end{lemma}

This immediately implies the following duality 
that was proved in \cite{Broadie85} by a different 
(and longer) argument. 

\begin{cor}[\cite{Broadie85}, Lemma 4.4] \label{pcduality}
The dual of a perfectly centered polytope is perfectly centered.
\end{cor}

The following theorem is equivalent to a theorem due to Broadie.

\begin{lemma}[\cite{Broadie85}, Theorem 2.1]\label{broadie}
If $P$ is a perfectly centered polytope, then $H$
is a facet of the Minkowski
sum $P+P^*$ if and only if $H$ is the sum of a face $F$ of $P$ with its
associated dual face $F^D$ in $P^*$.
\end{lemma}


Our first goal is to extend the characterization of facets
to all faces and to determine the face lattice of the 
Nesterov rounding $P+P^*$ of a perfectly centered polytope.

\begin{lemma}[\cite{Fukuda04}, Proposition 2.1]
Let $P_1,\ldots,P_k$ be polytopes in $\RRd$ and let $P=P_1+\cdots+P_k$.
Then a nonempty subset $F$ of $P$ is a face of $P$ if and only if $F=F_1+\cdots+F_k$
for some faces $F_i$ of $P_i$ such that there exists $c\in\RRd$ (not depending
on $i$) with $F_i=S(P_i;c)$ for all $i$. Furthermore, the decomposition
$F=F_1+\cdots+F_k$ of any nonempty face $F$ is unique.
\end{lemma}

\begin{lemma}\label{lem:subface}
Let $P$ be a perfectly centered polytope.
If a facet of $P+P^*$ is decomposed into two faces $F\subseteq P$
and $F^D \subseteq P^*$, then any nonempty subface $G$ of $F$ generates with
$F^D$ a subface of $F+F^D$ of dimension $\dim(G)+\dim(F^D)$.
\begin{proof}
This is the case because the faces $G$
and $F^D$ span affine spaces which are orthogonal to each other.
\end{proof}
\end{lemma}

In other words, for any two faces $F$ and $G$ of $P$ with $G\subseteq F$,
$G$ and $F^D$ sum to a face of $P+P^*$. We will show that there are no
other faces in $P+P^*$.

\begin{lemma}\label{lem:relat}
Let $P$ be a polytope. Let two nonempty faces of its Nesterov
rounding $P+P^*$ be decomposed as $G_1+F_1^D$ and $G_2+F_2^D$.
Then
$$
G_1+F_1^D\subseteq G_2+F_2^D \quad \Leftrightarrow \quad
G_1\subseteq G_2,\;F_1\supseteq F_2.
$$
\begin{proof}
Let two nonempty faces of its Nesterov
rounding $P+P^*$ be decomposed as $G_1+F_1^D$ and $G_2+F_2^D$.

If $G_1\subseteq G_2$ and $F_1 \supseteq F_2$, we have 
$F_1^D\subseteq F_2^D$, and thus $G_1+F_1^D\subseteq G_2+F_2^D$.

For the converse direction, observe that for two faces $A$ and $B$
of a polytope $P$,
$A\nsubseteq B$ if and only if $\cl(\NN(A;P))\cap\NN(B;P)=\emptyset$.
Assume $G_1\nsubseteq G_2$, that is, $\cl(\NN(G_1;P))\cap\NN(G_2;P)=\emptyset$.  This implies
\begin{align*}
& \cl(\NN(G_1+F_1^D;P+P^*))\cap\NN(G_2+F_2^D;P+P^*)\\
& = \cl(\NN(G_1;P)\cap\NN(F_1^D;P^*))\cap\NN(G_2;P)\cap\NN(F_2^D;P^*)\\
& \subseteq \cl(\NN(G_1;P))\cap\NN(G_2;P)\cap\cl(\NN(F_1^D;P^*))\cap\NN(F_2^D;P^*)=
\emptyset.
\end{align*}
Consequently, $G_1+F_1^D\nsubseteq G_2+F_2^D$. The same holds if
$F_1\nsupseteq F_2$ by symmetry.
\end{proof}
\end{lemma}

Now we are ready to prove:
\begin{quote}
\textbf{Theorem~\ref{mainthm}}~
Let $P$ be a perfectly centered polytope.
A subset $H$ of $P+P^*$ is a nontrivial face of  $P+P^*$ if 
and only if  $H=G+F^D$ for some ordered nontrivial faces 
$G\subseteq F$ of $P$.
\end{quote}
\begin{proof}
By Lemma~\ref{broadie}, the facets of $P+P^*$ are of form
$F+F^D$ for some nontrivial face $F$ of $P$.
Lemma~\ref{lem:subface} says that
 if $F$ and $G$ are nontrivial faces of $P$
with $G\subseteq F$, then $G+F^D$ is a face of the sum polytope.
Finally, Lemma~\ref{lem:relat} shows that all the faces are of that kind, since it proves
that there are no other subfaces to the facets.
\end{proof}

\begin{cor}\label{lattice}
The face lattice of the Nesterov
rounding $P+P^*$ of a perfectly centered polytope
is determined by that of $P$.
\end{cor}


\begin{thm} \label{thm:pcpc}
The Nesterov rounding of a perfectly centered polytope is also perfectly
centered.
\begin{proof}
Le $P$ be a perfectly centered polytope.
Let $F$ and $G$ be nontrivial faces of $P$ with $G\subseteq F$.
We'll denote by $m_F$ and $m_G$ the unique points in their intersections
with their respective normal cones.
By Theorem~\ref{mainthm}, it suffices to show that $m_G+m_{F^D}\in \NN(G;P)\cap\NN(F^D;P^*)$.
By Lemma~\ref{lem:coneinter}, $m_G\in\NN(G;P)\cap\NN(G^D;P^*)$. Also,
$m_{F^D}\in\NN(F;P)\cap\NN(F^D;P^*)$. Since $G\subseteq F$,
$\NN(F;P) \subseteq \cl(\NN(G;P))$. Since $m_G\in\NN(G;P)$
and $m_{F^D}\in \cl(\NN(G;P))$, $m_G+m_{F^D} \in\NN(G;P)$.
By symmetry, $m_G+m_{F^D} \in\NN(F^D;P^*)$, completing the proof.
\end{proof}
\end{thm}

\begin{figure}[htb]
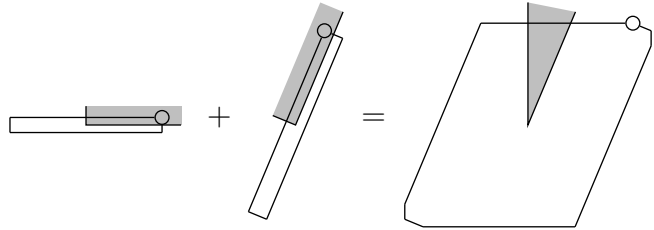

$$
\psset{xunit=0.5cm,yunit=0.5cm,shortput=nab,linewidth=0.5pt,arrowsize=2pt 3,labelsep=2pt}
\pspicture[.5](-3,-3)(3,3)
\psline[fillstyle=solid,fillcolor=lightgray,linetype=1,linestyle=none](2.5,0)(0,0)(0,0.5)(2.5,0.5)
\pnode(-2,-0.2){A}
\pnode(-2,0.2){B}
\cnode(2,0.2){1mm}{C}
\pnode(2,-0.2){D}
\ncline{A}{B}
\ncline{B}{C}
\ncline{C}{D}
\ncline{D}{A}
\psline(2.5,0)(0,0)(0,0.5)
\endpspicture
+
\pspicture[.5](-1.5,-3)(1.5,3)
\psline[fillstyle=solid,fillcolor=lightgray,linetype=1,linestyle=none](1.25,3)(0,0)(-0.6,0.25)(0.65,3.25)
\pnode(-0.76,-2.5){Ax}
\pnode(-1.24,-2.3){Bx}
\cnode(0.76,2.5){1mm}{Cx}
\pnode(1.24,2.3){Dx}
\ncline{Ax}{Bx}
\ncline{Bx}{Cx}
\ncline{Cx}{Dx}
\ncline{Dx}{Ax}
\psline(1.25,3)(0,0)(-0.6,0.25)
\endpspicture
=
\psset{xunit=0.5cm,yunit=0.5cm,shortput=nab,linewidth=0.5pt,arrowsize=2pt 3,labelsep=2pt}
\pspicture[.5](-3.5,-3)(3.5,3)
\psline[fillstyle=solid,fillcolor=lightgray,linetype=1,linestyle=none](1.25,3)(0,0)(0,3.25)
\pnode(-2.76,-2.7){Ay}
\pnode(-3.24,-2.5){By}
\pnode(-3.24,-2.1){Cy}
\pnode(-1.24,2.7){Dy}
\cnode(2.76,2.7){1mm}{Ey}
\pnode(3.24,2.5){Fy}
\pnode(3.24,2.1){Gy}
\pnode(1.24,-2.7){Hy}
\ncline{Ay}{By}
\ncline{By}{Cy}
\ncline{Cy}{Dy}
\ncline{Dy}{Ey}
\ncline{Ey}{Fy}
\ncline{Fy}{Gy}
\ncline{Gy}{Hy}
\ncline{Hy}{Ay}
\psline(1.25,3)(0,0)(0,3.25)
\endpspicture
$$
\caption{\label{fig:imperfsum}A non-perfectly-centered sum of perfectly centered polytopes}
\end{figure}

\begin{note}
The sum of two perfectly centered polytopes is not always perfectly centered.
For example, in Figure~\ref{fig:imperfsum}, both rectangles are perfectly centered, but their sum is not,
since the sum of the two marked vertices is not in its normal cone.

\end{note}

\subsection{The f-vector of the Nesterov rounding of a simplex}
Here we apply Theorem~\ref{mainthm} to perfectly centered simplices.

\begin{thm} \label{thm:simplexrounding}
Let $\Delta_d$ be a perfectly centered simplex of dimension $d$. Then,
the $f$-vector of the Nesterov rounding of $\Delta_d$ is given by
$$
f_k(\Delta_d+\Delta_d^*)=
\left(\begin{array}{c}
  d+1 \\
  k+2
\end{array}\right)
\left(
2^{k+2}-2
\right), \quad \text{ for } 0 \le k \le d-1.
$$
\begin{proof}
Let $\Delta_d$ be a perfectly centered simplex of dimension $d$.
The $f$-vector of $\Delta_d$ is given by
$$
f_k(\Delta_d)=\left(\begin{array}{c}d+1 \\ k+1\end{array}\right), 
\quad \text{ for } 0 \le k \le d-1.
$$
By Theorem~\ref{mainthm}, the faces of $\Delta_d+\Delta_d^*$
can be characterized as the sums $F^D+G$, with $G\subseteq F$ non-trivial
faces of $\Delta_d$.

Let $S$ and $T$ be the vertex sets of respectively $G$ and $F$, with
$S\subseteq T$, and denote $U=T\setminus S$. The dimension $k$ of
$F^D+G$ is $\dim(F^D)+\dim(G)= d-1+\dim(G)-\dim(F)=d-1+|S|-|T|=d-1-|U|$.

So the number of faces of dimension $k$ can be written as $pq$, where $p$
is the number of possible choices
of $U$ with $|U|=d-1-k$, and $q$ is the number of choices of $S$ non-empty,
so that $S\cap U = \emptyset$ and $|T|=|S\cup U| < d+1$. Thus we have
$$
p=\left(\begin{array}{c}
  d+1 \\
  k+2
\end{array}\right)
\quad\quad\mbox{and}\quad\quad
q=2^{k+2}-2.
$$
\end{proof}
\end{thm}

\subsection{The f-vector of the Nesterov rounding of a cube}
\begin{thm} \label{thm:cuberounding}
Let $\Box_d$ be a cube of dimension $d$.  Then,
the $f$-vector of the Nesterov rounding of $\Box_d$ is given by
$$
f_k(\Box_d+\Box_d^*)=
\left(\begin{array}{c}
  d \\
  k+1
\end{array}\right)
2^{d-k-1}
\left(3^{k+1}-1\right), \quad \text{ for } 0 \le k \le d-1.
$$
\begin{proof}
Let $\Box_d$ be a cube of dimension $d$.
Then $\Box_d$ has $3^n-1$ nontrivial faces,
which can be decomposed as:
$$
f_k(\Box_d)=
\left(\begin{array}{c}
  n \\
  d
\end{array}\right)
2^{n-d}, \quad \text{ for } 0 \le k \le d-1.
$$
By Theorem~\ref{mainthm}, the faces of $\Box_d+\Box_d^*$
can be characterized as the sums $F^D+G$, with $G\subseteq F$ non-trivial
faces of $\Box_d$.

Let $S$ and $T$ be the sets of fixed coordinates of respectively
$G$ and $F$, with $T\subseteq S$, and denote $U=S\setminus T$.
The dimension $k$ of $F^D+G$ is
$\dim(F^D)+\dim(G)= d-1+\dim(G)-\dim(F)=d-1+(d-|S|)-(d-|T|)=d-1-|U|$.

So the number of faces of dimension $k$ can be written as $pqr$, where $p$
is the number of possible choices of $U$ with $|U|=d-1-k$,
$q$ is the number of ways to fix the coordinates in $U$, and
$r$ is the number of choices of $G$, so that $S\cap U = \emptyset$
and $|T|=|S\cup U| < d+1$. We have
$$
p=\left(\begin{array}{c}
  d \\
  k+1
\end{array}\right)
,\quad\quad
q=2^{d-k-1}
\quad\quad
\mbox{and}
\quad\quad
r=3^{k+1}-1.
$$
\end{proof}
\end{thm}

\subsection{Repeated Nesterov rounding in dimension 3}
We'll use the following notation: $f_k^{(i)}(P)$ denotes the number of
$k$-dimensional faces in a polytope $P$ after executing the Nesterov rounding
$i$ times.

\begin{thm} \label{thm:repeatedNR}
Let $P$ be a perfectly centered $3$-dimensional polytope $P$.
Then the following relations hold:
$$
f_0^{(n)} = 4^{n-1} f_0^{(1)},
$$
$$
f_1^{(n)} = 2\cdot4^{n-1} f_0^{(1)},\;\mbox{and}
$$
$$
f_2^{(n)} = f_2^{(1)}+(4^{n-1}-1) f_0^{(1)}.
$$
\begin{proof}
Let $P$ be a perfectly centered three-dimensional polytope.
By Corollary~\ref{lattice},
$$
f_2^{(n)} = f_0^{(n-1)}+f_1^{(n-1)}+f_2^{(n-1)}.
$$
It is a general property of face lattices that for two faces $G\subseteq F$
so that $\dim(G)+2=\dim(F)$ there are exactly two faces $H_1$ and $H_2$ of
dimension $\dim(G)+1$ so that $G\subseteq H_1\subseteq F$ and
$G\subseteq H_2\subseteq F$.
In a Nesterov rounding, it means that all $(d-3)$-dimensional faces,
which are sums of a face $G$ and $F^D$, $G\subseteq F$
so that $\dim(G)+2=\dim(F)$ are contained
in four $(d-2)$-dimensional faces, which are $G+H_1^D$, $G+H_2^D$, $H_1+F^D$
and $H_2+F^D$, and four $(d-1)$-dimensional faces, which are $G+G^D$,
$H_1+H_1^D$, $H_2+H_2^D$ and $F+F^D$.
In the $3$-dimensional case, it means that all vertices are contained
in four incident edges and four facets. Since each edge contains exactly
two vertices, we have:
$$
f_1^{(n)}=2f_0^{(n)},\;\forall n\geq 1.
$$
Since the number of vertices in the next Nesterov rounding is equal to the
number of pairs of a vertex and its containing facets, it also means that:
$$
f_0^{(n+1)}=4f_0^{(n)},\;\forall n\geq 1.
$$
Thus we have the following equations:
$$
f_0^{(n)} = 4^{n-1} f_0^{(1)},
$$
$$
f_1^{(n)} = 2\cdot4^{n-1} f_0^{(1)},\;\mbox{and}
$$
$$
f_2^{(n)} = f_2^{(n-1)}+3\cdot 4^{n-2} f_0^{(1)}
$$
$$
\Rightarrow f_2^{(n)} = f_2^{(1)}+(4^{n-1}-1) f_0^{(1)}.
$$
\end{proof}
\end{thm}

Note that the ratio of the number of facets to that of vertices tends towards $1$.

\section{Maximizing faces}
It is natural to explore possible bounds for the number of faces in
Minkowski sums of polytopes. The description of a Minkowski sum can
be exponential in terms of the description (binary) size of the summands. For
instance, the sum of $d$ orthogonal segments in $d$-dimensional space
is the $d$-hypercube, which has $2^d$ vertices, but only $2d$ facets.

In this section, we obtain some tight bounds on the number of faces in Minkowski
sums, in terms of number of vertices in the summands, and of the dimension.

\subsection{Bounds on vertices}
We will show an upper bound for the number of vertices in a Minkowski sum,
then we will show this bound is attainable provided the dimension is big enough
in relation to the number of polytopes.

Each vertex in a Minkowski sum is decomposed into a sum of vertices
of the summands. Since each vertex has a different decomposition,
we arrive to the following trivial upper bound:

\begin{lemma}[Trivial upper bound] \label{vertices}
Let $P_1,\ldots,P_n$ be polytopes. Then the following gives 
an upper bound on the number of vertices of their Minkowski sum:
$$
f_0(P_1+\cdots+P_n) \leq \prod_{i=1}^n f_0(P_i).
$$
\end{lemma}

Now we are ready to prove:
\begin{quote}
\textbf{Theorem~\ref{thmvertices}}~
In dimension $d\geq 3$, it is possible to choose
$(d-1)$ polytopes so that the trivial upper bound for vertices is attained.
\end{quote}
\begin{proof}
Let $P_i$, $i=1,\ldots,d-1$, be $d$-dimensional polytopes, and $v_{i,j}$
their vertices, $j=1,\ldots,n_i$ where $n_i\geq 1$ is the number
of vertices of the polytope $P_i$. We set the coordinates of the
vertices to be:
$$
v_{i,j} = \cos\left(\frac{j}{n_i+1}\pi\right)\cdot \boldsymbol{e_i} +
          \sin\left(\frac{j}{n_i+1}\pi\right)\cdot \boldsymbol{e_d}.
$$
Where $\boldsymbol{e_j}$'s are the unit vectors of an orthonormal basis
of the $d$-dimensional space.
So the vertices of $P_i$ are placed on the unit half-circle in
the space generated by $\boldsymbol{e_i}$ and $\boldsymbol{e_d}$.
Observe that the polytopes are two-dimensional for now. 
By the construction, one can easily verify
that
$$
v_{i,j}\in \NN(\{v_{i,j}\};P_i).
$$
This stays true if we add anything to those vectors in the spaces orthogonal to that of the half-circle:
$$
v_{i,j}+\sum_{k\neq i,d}\alpha_k\boldsymbol{e_k}\in \NN(\{v_{i,j}\};P_i),
\forall \alpha_k\in \mathbb{R}.
$$
So for any choice of $S=\{j_i\}_{i=1}^{d-1}$, $j_i=1,\ldots,n_i$,
we can build this vector:
$$
v_S = \sum_{i=1}^{d-1}\cot\left(\frac{j_i}{n_i+1}\pi\right)\cdot
\boldsymbol{e_i} + \boldsymbol{e_d}.
$$
This vector $v_S$, projected
to the space generated by $\boldsymbol{e_d}$ and any $\boldsymbol{e_i}$, is equal to
$
\cot\left(\frac{j_i}{n_i+1}\pi\right)\cdot \boldsymbol{e_i} + \boldsymbol{e_d}
$ which is collinear with $\cos\left(\frac{j_i}{n_i+1}\pi\right)\cdot \boldsymbol{e_i} +
\sin\left(\frac{j_i}{n_i+1}\pi\right)\cdot \boldsymbol{e_d}$, and thus belongs
to $\NN(\{v_{i,j_i}\};P_i)$. So we have that:
$$
v_S \in \bigcap_{i=1}^{d-1}\NN(\{v_{i,j_i}\};P_i),
$$
and since this intersection is not empty, it means that
$v_{j_1},\ldots,v_{j_{d-1}}$ is a vertex of the Minkowski sum
$P_1+\cdots+P_{d-1}$. This stays true for any choice of
$S=\{j_i\}_{i=1}^{d-1}$, so the Minkowski sum has $\prod_{i=1}^{d-1}n_i$
vertices. The polytopes $P_i$ thus defined are 2-dimensional.
The property still stands if we add small perturbations to the vertices
to make the polytopes full-dimensional.
\end{proof}

\subsection{Bounds on facets}
It appears to be much harder to find tight upper bounds on facets of Minkowski
sums than on vertices. This is due to the fact that vertices of the sum
decompose only in sum of vertices of the summands, while facets decompose in
faces that can have any dimension. Results are therefore limited for now
to low-dimensional cases.

Let us recall now Theorem~\ref{thmfaces}, which gives the number
of facets in a 3-dimensional Minkowski sum of polytopes relatively in
general position:
\begin{quote}
\textbf{Theorem~\ref{thmfaces}}~
Let $P_1,\ldots,P_n$ be 3-dimensional polytopes relatively in general
position, and $P$ their sum. Then the following equations hold:
$$
2 f_2(P)-f_1(P) = \sum_{i=1}^n (2 f_2(P_i)- f_1(P_i)) ,
$$
$$
f_2(P)-f_0(P)+2 = \sum_{i=1}^n (f_2(P_i)-  f_0(P_i)+2) ,
$$
$$
f_1(P)-2f_0(P)+4 = \sum_{i=1}^n (f_1(P_i)-2 f_0(P_i)+4 ).
$$
\end{quote}
\begin{proof}
Let $P_1,\ldots,P_n$ be 3-dimensional polytopes relatively in general
position and $P$ their sum.
A facet of $P$ can be either \emph{pure}, which means its decomposition
contains exactly one facet of one of the summands and vertices otherwise,
or it can be \emph{mixed}, which means the decomposition contains exactly
two edges and vertices otherwise.

In terms of normal cones, the normal
ray of a pure facet is the intersection of the normal ray of a single facet
with $3$-dimensional normal cones of vertices. This means each facet in
the summands will generate exactly one pure facet of the sum:
$$
f_2^{pure}(P) = f_2(P_1)+\cdots+f_2(P_n).
$$

The normal ray of a mixed facet is the intersection of the
$2$-dimensional normal cones of exactly two edges with $3$-dimensional
normal cones of vertices.
This means that each of those two normal cones are split into two
non-connected sets by removal of the intersection. This creates in
effect four normal cones of four different edges in the Minkowski sum.
So every occurence of a mixed facet augments by exactly two
the number of edges in the sum:
$$
2f_2^{mix}(P) = f_1(P)-(f_1(P_1)+\cdots+f_1(P_n)).
$$

Combining the two equations, we get the first part of the theorem.
The two other parts are deduced using Euler's equation:
$f_0+f_2=f_1+2.$
\end{proof}

$ $\\
This result allows us to find tight upper bounds on the number of
edges and facets when summing two $3$-dimensional polytopes:

\begin{cor} \label{cor:facets}
Let $P_1$ and $P_2$ be polytopes in dimension 3, and $P=P_1+P_2$ their
Minkowski sum. Then we have the following tight bounds:
$$
f_2(P)\leq f_0(P_1)f_0(P_2)+ f_0(P_1)+ f_0(P_2)-6,
$$
$$
f_1(P)\leq 2f_0(P_1)f_0(P_2)+ f_0(P_1)+ f_0(P_2)-8.
$$
\begin{proof}
We first need to prove that the maximum number of facets can be
attained if $P_1$ and $P_2$ are relatively in general position.

Let $P_1$ and $P_2$ be polytopes in dimension 3, so that their
Minkowski sum $P$ has the maximal number of facets. They
are not relatively in general position if and only if
a facet of the summands is summing with a facet or an edge.
Suppose we perturb $P_2$ by a small well-chosen rotation.
The facets previously contained in a non-exact decomposition will now sum
to a facet with an exact decomposition. Therefore, there will be at least
as many facets as before in the sum, and the summands will now be
relatively in general position.

We now find what the maximum is when summands are relatively
in general position. By Theorem~\ref{thmfaces}, it is sufficient
to maximize $f_0(P)$ and $f_2(P_i)-f_0(P_i)$ for each $i$. We can
do this by using simplicial polytopes disposed as indicated in
Theorem~\ref{thmvertices}. Since vertices and facets are maximized
in the sum, so are edges.
\end{proof}
\end{cor}

\subsection{Sums of cyclic polytopes}

We will show here an upper bound for the number of faces of each dimension
a Minkowski sum of polytopes can have. We will then show this bound is attained for
lower dimensions by certain sums of cyclic polytopes.

\begin{lemma}[Trivial upper bound 2] \label{maxfaces}
Let $P_1,\ldots,P_n$ be $d$-dimensional polytopes. For each $k=0,\ldots,d-1$
and $n\geq 1$, the number of $k$-faces of $P_1+\cdots+P_n$ is bounded by:
$$
f_k(P_1+\cdots+P_n) \leq
\sum_{
\begin{array}{c}
1\leq s_i \leq f_0(P_i)\\
s_1+\cdots+s_n=k+n
\end{array}}
\prod_{i=1}^n
\left(\begin{array}{c}
f_0(P_i)\\
s_i
\end{array}\right),
$$
where $s_i$'s are integral.
\begin{proof}
Let $P_1,\ldots,P_n$ be $d$-dimensional polytopes,
and $F$ a $k$-dimensional face of $P_1+\cdots+P_n$.
Let $F_i\subseteq  P_i,\;i=1,\ldots,n$ be the decomposition of $F$.
Let $k_1,\ldots,k_n$ be the dimensions of respectively $F_1,\ldots,F_n$.
Then $k_1+\cdots+k_n\geq k$. The minimal number of vertices for a face
of dimension $k_i$ is $k_i+1$. So the total number of vertices contained
in faces of the decomposition of $F$ is at least $k+n$. For any fixed
$k_1,\ldots,k_n$, the number of possible choices of $s_i$ vertices
for each $P_i$ is:
$$
\prod_{i=1}^n
\left(\begin{array}{c}
f_0(P_i)\\
s_i
\end{array}\right).
$$
\end{proof}
\end{lemma}
\noindent
Note that the above lemma is an extension of Lemma~\ref{vertices}.

Cyclic polytopes are known to have the maximal number of faces for
any fixed number of vertices. This property is somewhat carried
on to their Minkowski sum.
We define the \emph{moment curve} as the curve in the $d$-dimensional
space which is the set of points of form $(x,x^2,x^3,\ldots,x^d)$.
We call $P$ a \emph{cyclic polytope} if its vertices are all on the
moment curve. Cyclic polytopes have the following properties: Their number
of faces is maximal, for faces of all dimensions, over all polytopes with this
dimension and this number of vertices \cite{mcmullen70}. They are also
\emph{simplicial}, i.e. all their faces are simplices.
Moreover, they are \emph{$\lfloor\frac{d}{2}\rfloor$-neighbourly},
which means that the convex hull of any set of
$\lfloor\frac{d}{2}\rfloor$ vertices
of a cyclic polytope $P$ is a face of $P$ \cite[4.7]{Grunbaum67}.
Since each face is a simplex, this also means that any set of
$\lfloor\frac{d}{2}\rfloor$ vertices are affinely independent.
For more details concerning cyclic polytopes, please refer to
\cite{Grunbaum67}.

Note that if we choose a set $S$ of points on the moment curve,
with $|S|\leq\lfloor\frac{d}{2}\rfloor$, $\conv(S)$ will form a face
of any polytope $P$ having $S$ as a subset of its vertices,
no matter how the other vertices are chosen.
That is, there is always a linear function $\langle m_S,x\rangle$ so that
$S(P;m_S) = \conv(S)$.

Let us now recall:
\begin{quote}
\textbf{Theorem~\ref{thmcyc}}~
In dimension $d\geq 4$,
it is possible to choose $n\leq\lfloor\frac{d}{2}\rfloor$
polytopes $P_1$, $P_2$, $\ldots$, $P_n$ so that the trivial upper bound for 
the number of $k$-faces
of $P_1+ \cdots + P_n$ is attained for all $0\le k \le \lfloor\frac{d}{2}\rfloor-n$.
\end{quote}
\begin{proof}
Let $P$ be the Minkowski sum of polytopes $P_1,\ldots,P_n$
whose vertices are all distinct on the moment curve, with
$k=\lfloor\frac{d}{2}\rfloor-n$, $k\in \mathbb{N}$.

Let $S_1\subseteq\mathcal{V}(P_1),\ldots,S_n\subseteq\mathcal{V}(P_n)$
be subsets of the vertices of the polytopes such that
$S_i\neq\emptyset,\;\forall i$ and $|S_1|+\cdots+|S_n|=k+n$.
Since $k+n\leq\lfloor\frac{d}{2}\rfloor$, there is a linear function
maximized at $S_1,\ldots,S_n$ on the moment curve. Therefore,
$\conv(S_i)$ is an $(|S_i|-1)$-dimensional face of $P_i$, $\forall i=1,\ldots,n$.
Since the same linear function is maximized over each $P_i$ on these faces,
they sum up to a face of $P$. Since the set of vertices $S_1\cup\ldots\cup S_n$
is affinely independant, $\dim(\conv(S_1)+\cdots+\conv(S_n)) = \dim(\conv(S_1))+\cdots+\dim(\conv(S_n)) = |S_1|+\cdots+|S_n|-n=n+k-n = k$.
\end{proof}

\begin{ack}
We would like to thank one of the referees for bringing Broadie's paper
\cite{Broadie85} to our attention. We are also grateful to G\"unther Ziegler
who provided us with much simpler proofs for Theorems
\ref{thm:simplexrounding} and \ref{thm:cuberounding}.
\end{ack}

\bibliographystyle{abbrv}
\bibliography{biblio}

\end{document}